%% file: IRPA2.tex
\documentclass[letterpaper,11pt,twoside,keywordsasfootnote,addressatend,noinfoline]{article}

\usepackage{imsart}
\usepackage{fullpage}
\usepackage[english]{babel}
\usepackage{amssymb}
\usepackage{amsmath}
\usepackage{amsthm}
\usepackage{epsfig}
\usepackage{subfigure}
\usepackage{mathrsfs}
\usepackage{color}

\newtheorem{theorem}{Theorem}
\newtheorem{lemma}[theorem]{Lemma}

\renewenvironment{proof}{\noindent{\bf Proof.}}{\hspace*{2mm}~$\square$}
\newenvironment{proofof}[1]{\noindent{\bf Proof of #1.}}{\hspace*{2mm}~$\square$}

\newcommand{\B}{\mathscr{B}}

\newcommand{\Z}{\mathbb{Z}}

\newcommand{\lat}{\mathscr{L}}

\newcommand{\ind}{\mathbf{1}}
\newcommand{\ep}{\epsilon}
\newcommand{\n}{\hspace*{-5pt}}
\newcommand{\bul}[1]{#1\,-\,$\bullet$}

\newcommand{\arrow}[1]{\overset{\hbox{\tiny #1}}{\longrightarrow}}

\DeclareMathOperator{\binomial}{Binomial}
\DeclareMathOperator{\geometric}{Geometric}
\DeclareMathOperator{\poisson}{Poisson}
\DeclareMathOperator{\exponential}{Exponential}
\DeclareMathOperator{\card}{card}

\begin{document}

\begin{frontmatter}
\title{Multiple transition points for an epidemic model}
\runtitle{Multiple transition points for an epidemic model}
\author{Lamia Belhadji and Nicolas Lanchier}
\runauthor{Lamia Belhadji and Nicolas Lanchier}
\address{Faculty of Exact Sciences \\ and Computer Science \\ University of Abdelhamid Ibn Badis \\ Mostaganem, Algeria \\ lamia.belhadji@univ-mosta.dz}
\address{School of Mathematical \\ and Statistical Sciences \\ Arizona State University \\ Tempe, AZ 85287, USA. \\ nicolas.lanchier@asu.edu}
\maketitle

\begin{abstract} \
 This paper is concerned with a spatial epidemic model on the $d$-dimensional square lattice that extends the contact process in two ways by assuming that each site is occupied by a social cluster of $N$ individuals rather than a single individual, and that infected individuals can be either asymptomatic or symptomatic. Newly infected individuals begin as asymptomatic and develop the symptoms at a spontaneous mutation rate $\gamma$, and all infected individuals recover at rate one. The model also distinguishes between external~(between adjacent clusters) and internal~(within clusters) infection rates for each state. Our analysis shows that, when either asymptomatic or symptomatic individuals only spread the disease externally and at a subcritical rate, while the other class only transmits the disease internally and at a large rate, the system exhibits at least two phase transitions with respect to the mutation rate~$\gamma$ when social clusters are large. Specifically, the epidemic dies out for extreme values but survives for intermediate values of the mutation rate.
\end{abstract}



\end{frontmatter}


\section{Introduction}
\label{sec:intro}
\begin{figure}[h!]
\centering
\scalebox{0.65}{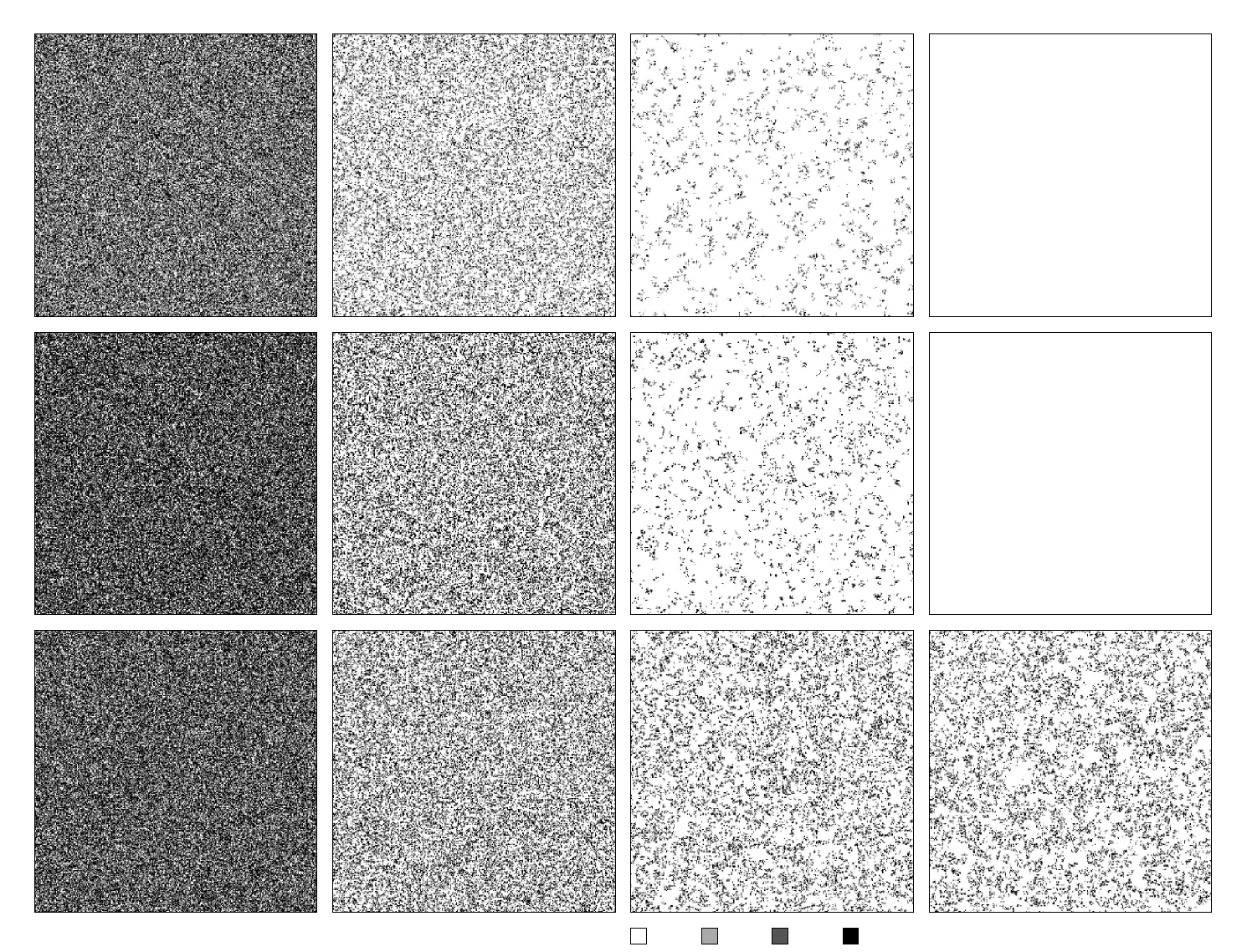}
\caption{\upshape
 Snapshots of the process on the~$300 \times 300$ torus at various times and for various values of~$\gamma$.
 In all three pictures, the parameters are~$\lambda_1 = 1/4 < \lambda_c, \lambda_2 = \phi_1 = 0$ and~$\phi_2 = N = 3$, and the process starts with one heathy, one asymptomatic, and one symptomatic individuals in each cluster.
 The two simulations at the top illustrate the first part of Theorem~\ref{th:gamma1}, while the simulation at the bottom illustrates the second part of the theorem, and the three simulations together show the presence of at least two phase transitions with respect to~$\gamma$.}
\label{fig:IPS}
\end{figure}

 This paper shows the lack of monotonicity, resulting in the existence of~(at least) two phase transitions, for the epidemic model recently introduced in~\cite{belatreche_belhadji_lanchier_2026}.
 This model is a variant of Harris' contact process~\cite{harris_1974} on the~$d$-dimensional square lattice that replaces each lattice point/single individual by a social cluster of~$N$ individuals, like the model in~\cite{belhadji_lanchier_2005}, and distinguishes infected individuals depending on whether they are asymptomatic or symptomatic, like the model in~\cite{belhadji_lanchier_mercer_2024,lanchier_2026}.
 For a review of these two epidemic models, we refer the reader to~\cite[Chapter~2]{lanchier_2024}.
 In particular, the state of the process at time~$t$ is a spatial configuration
 $$ \xi_t : \Z^d \to \{(i, j) \in \Z_+ \times \Z_+ : i + j \leq N \}, $$
 where~$\xi_t (x) = (i, j)$ means that social cluster~$x$ contains~$i$ asymptomatic individuals, and~$j$ symptomatic individuals, while the other individuals are healthy.
 The presence of social clusters motivates the introduction of two infection rates depending on whether the infection originates from a nearby cluster or from the same cluster.
 Similarly, the distinction between asymptomatic and symptomatic individuals motivates the introduction of two infection rates depending on the state of the individual the infection originates from.
 Combining the two gives four infection rates: external asymptomatic infections at rate~$\lambda_1$, external symptomatic infections at rate~$\lambda_2$, internal asymptomatic infections at rate~$\phi_1$, and internal symptomatic infections rate~$\phi_2$.
 In addition, we assume that newly infected individuals are always asymptomatic and become symptomatic at the mutation rate~$\gamma$, that infected individuals recover at rate one, and that external infections are neglected once a cluster is infected.
 In particular, writing~$y \sim x$ to indicate that cluster~$y$ is one of the~$2d$ nearest neighbors of cluster~$x$, the transition rates at cluster~$x$ are given by
 $$ \begin{array}{rclcl}
      (0, 0) & \n \to \n & (1, 0)         & \quad \hbox{at rate} \quad & \sum_{y \sim x} \sum_{i, j} \,(i \lambda_1 + j \lambda_2) \,\ind \{\xi_t (y) = (i, j) \}, \vspace*{4pt} \\
      (i, j) & \n \to \n & (i + 1, j)     & \quad \hbox{at rate} \quad & i \phi_1 + j \phi_2 \quad \hbox{for} \ i + j < N, \vspace*{4pt} \\
      (i, j) & \n \to \n & (i - 1, j + 1) & \quad \hbox{at rate} \quad & i \gamma \quad \hbox{for} \ i + j \leq N, j < N, \vspace*{4pt} \\
      (i, j) & \n \to \n & (i - 1, j)     & \quad \hbox{at rate} \quad & i \quad \hbox{for} \ i + j \leq N, \vspace*{4pt} \\
      (i, j) & \n \to \n & (i, j - 1)     & \quad \hbox{at rate} \quad & j \quad \hbox{for} \ i + j \leq N. \end{array} $$
 From top to bottom, these transitions are interpreted as external infection of a healthy cluster, additional internal infection of an infected cluster, mutation from asymptomatic to symptomatic, recovery of an asymptomatic individual, and recovery of a symptomatic individual.
 We assume that the epidemic model starts from a translation-invariant distribution with a positive density of asymptomatic individuals and a positive density of symptomatic individuals.
 Like the contact process, there are two possible long-term behaviors:
 $$ \begin{array}{rcrcl}
      \hbox{survival} & \hbox{when} & \liminf_{t \to \infty} P (\xi_t (x) \neq (0, 0)) > 0 & \hbox{for all} & x \in \Z^d, \vspace*{4pt} \\
    \hbox{extinction} & \hbox{when} &    \lim_{t \to \infty} P (\xi_t (x) \neq (0, 0)) = 0 & \hbox{for all} & x \in \Z^d. \end{array} $$
 Our previous work~\cite{belatreche_belhadji_lanchier_2026} shows the existence of phase transitions between extinction and survival with respect to each of the four infection rates and the cluster size.
 Even though coupling methods comparing processes with different parameters fail to prove monotonicity, the analysis in~\cite{belatreche_belhadji_lanchier_2026} suggests monotonicity in the sense that, if there is survival for a given set of parameters, increasing the infection rates and the cluster size still results in survival, which is further supported by numerical simulations of the two-dimensional process.
 In contrast, increasing the mutation rate~$\gamma$ can either promote survival or promote extinction depending on the other parameters.
 More interestingly, based on the analysis of the mean-field model, we conjectured that, for some values of the infection rates and the cluster size, there is extinction for~$\gamma$ small or large, but survival for intermediate values, thus resulting in the existence of multiple phase transitions with respect to the mutation rate.
 The objective of this paper is to prove this conjecture. \vspace*{8pt} \\
\noindent
{\bf Mean-field model}.
 To hypothesize how to choose the infection parameters, we first return to the mean-field model, which consists of a system of coupled differential equations for the densities~$u_{ij}$ of clusters with~$i$ asymptomatic and~$j$ symptomatic individuals, and is obtained by assuming that the population is homogeneously mixing.
 Defining survival in the context of the mean-field model as the existence of a fixed point with a positive density of infected clusters,~\cite[Proposition~1]{belatreche_belhadji_lanchier_2026} shows that, for clusters of size two, there is survival if and only if
\begin{equation}
\label{eq:MF0}
\begin{array}{rcl}
  (\lambda_1 + \gamma \lambda_2)((1 + \phi_1)(1 + \phi_2) + (1 + \phi_1) + \gamma (1 + \phi_2)) \vspace*{4pt} \\ \hspace*{120pt} -
  (\gamma \phi_1 + \phi_2 + (1 + \gamma)(2 + \gamma)) & \n > \n & 0.
\end{array}
\end{equation}
 Using this preliminary result, we also proved that the mean-field model is monotone with respect to the infection rates in the sense that, if there is survival for a given set of parameters, increasing the infection rates still results in survival, and by contraposition, if there is extinction for a given set of parameters, decreasing the infection rates still results in extinction.
 The mean-field model, however, is not monotone with respect to the mutation rate.
 Indeed, taking~$\lambda_2 = \phi_1 = 0$, the condition for survival~\eqref{eq:MF0} reduces to
 $$ - \gamma^2 + (\lambda_1 (1 + \phi_2) - 3) \gamma + (\lambda_1 - 1)(2 + \phi_2) > 0. $$
 Then, taking for instance~$\lambda_1 = 3/4$ and~$\phi_2 = 11$, the left-hand side becomes
\begin{equation}
\label{eq:MF1}
    - \gamma^2 + 6 \gamma - 13/4 = - (\gamma - \gamma_-)(\gamma - \gamma_+) \quad \hbox{where} \quad \gamma_{\pm} = 3 \pm \sqrt{23}/2 > 0,
\end{equation}
 showing the existence of two phase transitions with
 $$ \hbox{extinction if} \ \gamma \leq \gamma_-, \quad \hbox{survival if} \ \gamma_- < \gamma < \gamma_+, \quad \hbox{extinction if} \ \gamma \geq \gamma_+. $$
 Similarly, when~$\lambda_1 = \phi_2 = 0$, the condition for survival~\eqref{eq:MF0} reduces to
 $$ (\lambda_2 - 1) \gamma^2 + (2 \lambda_2 (1 + \phi_1) - \phi_1 - 3) \gamma - 2 > 0. $$
 Taking~$\lambda_2 = 3/4$ and~$\phi_1 = 11$, the left-hand side becomes
\begin{equation}
\label{eq:MF2}
    - \gamma^2 / 4 + 4 \gamma - 2 = - (\gamma - \gamma_-')(\gamma - \gamma_+') / 4 \quad \hbox{where} \quad \gamma_{\pm}' = 8 \pm 2 \sqrt{14} > 0,
\end{equation}
 showing again the existence of two phase transitions with
 $$ \hbox{extinction if} \ \gamma \leq \gamma_-', \quad \hbox{survival if} \ \gamma_-' < \gamma < \gamma_+', \quad \hbox{extinction if} \ \gamma \geq \gamma_+'. $$
 In conclusion, at least in the mean-field model with clusters of size two, we have two phase transitions with respect to the mutation rate when asymptomatic individuals infect nearby clusters occasionally but not their own cluster, while symptomatic individuals cannot infect nearby clusters but are highly infectious in their own cluster, leading to scenario~\eqref{eq:MF1}.
 Similarly, there are also two phase transitions when asymptomatic individuals cannot infect nearby clusters but are highly infectious in their own cluster, while symptomatic individuals infect nearby clusters occasionally but not their own cluster, leading to scenario~\eqref{eq:MF2}. \vspace*{8pt} \\
\noindent
{\bf Stochastic spatial model}.
 Following the results for the mean-field model, to mimic scenario~\eqref{eq:MF1} above but in the context of the interacting particle system, we again take~$\lambda_2 = \phi_1 = 0$.
 To ensure extinction for~$\gamma$ small in the mean-field model, we need~$\lambda_1 < 1$.
 For the interacting particle system, this condition can be relaxed taking~$\lambda_1$ smaller than the critical value~$\lambda_c$ of the basic contact process.
 Finally, instead of simply taking~$\phi_2$ large, we set~$\phi_2 = N$ and study the process when the cluster size~$N$ is large.
 This choice of the four infection rates and the cluster size leads to the existence of two phase transitions with respect to~$\gamma$ when~$N$ is sufficiently large.
\begin{theorem}
\label{th:gamma1}
 Assuming that the infection rates are given by
 $$ \lambda_1 \in (0, \lambda_c), \quad \lambda_2 = 0, \quad \phi_1 = 0, \quad \phi_2 = N = \hbox{cluster size}, $$
 there are at least two phase transitions with respect to~$\gamma$ in the sense that
 \begin{itemize}
\item[a.]
 $\forall \,N > 1, \,\exists \,0 < \gamma_- < \gamma_+ < \infty : \gamma \not \in (\gamma_-, \gamma_+) \ \Rightarrow$ extinction, \vspace*{4pt}
\item[b.]
 $\forall \,\gamma \in (0, \infty), \,\exists \,N_+ < \infty : N > N_+ \ \Rightarrow$ survival.
\end{itemize}
\end{theorem}
\noindent
 Note that the theorem indeed implies the presence of at least two phase transitions for some range of the infection parameters and cluster size, with extinction for extreme values and survival for intermediate values of~$\gamma$.
 To see this, fix~$\gamma_0 \in (0, \infty)$, and apply the second part of the theorem to get an~$N_0$ finite such that the process survives, then apply the first part of the theorem to get two values~$\gamma_- < \gamma_0 < \gamma_+$ such that the process dies out for~$\gamma < \gamma_-$ and~$\gamma > \gamma_+$.
 This shows that, when the cluster size is~$N_0$, there is a phase transition from extinction to survival between~$\gamma_-$ and~$\gamma_0$, and another phase transition from survival to extinction between~$\gamma_0$ and~$\gamma_+$.
 Figure~\ref{fig:IPS} illustrates this result for clusters with only three individuals.
 Scenario~\eqref{eq:MF2} above can also be reproduced for the interacting particle system by simply reversing the roles of the asymptomatic and symptomatic infection rates, which gives the following theorem.
\begin{theorem}
\label{th:gamma2}
 Assuming that the infection rates are given by
 $$ \lambda_1 = 0, \quad \lambda_2 \in (0, \lambda_c), \quad \phi_1 = N = \hbox{cluster size}, \quad \phi_2 = 0, $$
 there are at least two phase transitions with respect to~$\gamma$ in the sense that
 \begin{itemize}
\item[a.]
 $\forall \,N > 1, \,\exists \,0 < \gamma_- < \gamma_+ < \infty : \gamma \not \in (\gamma_-, \gamma_+) \ \Rightarrow$ extinction, \vspace*{4pt}
\item[b.]
 $\forall \,\gamma \in (0, \infty), \,\exists \,N_+ < \infty : N > N_+ \ \Rightarrow$ survival.
\end{itemize}
\end{theorem}
\noindent
 Even though both theorems give the same conclusion, with very similar assumptions, the proofs are significantly different.
 This is due to the fact that newly infected individuals are always asymptomatic, while symptomatic individuals can only appear later as the result of a mutation, which completely breaks down symmetry.
 In particular, reversing the values of the asymptomatic and symptomatic infection rates does not result in general in the same behavior.
 In fact, survival for intermediate values of~$\gamma$ follows from the quick spread of the disease across space in Theorem~\ref{th:gamma1}, but from the long persistence of the disease within clusters in Theorem~\ref{th:gamma2}.


\section{Graphical representation}
\label{sec:GR}
 The first step to prove the theorems is to construct the process graphically.
 As pointed out in~\cite{belatreche_belhadji_lanchier_2026}, it is convenient to view each social cluster as a pile of~$N$ sites corresponding to the individuals living in that cluster, and to think of the process as keeping track of the state of each individual, rather than the state of each cluster.
 In particular, the state at time~$t$ becomes
 $$ \xi_t : \Delta \longrightarrow \{0, 1, 2 \} \quad \hbox{where} \quad \Delta = \Z^d \times \{1, 2, \ldots, N \}, $$
 and where~0 means healthy, 1 means asymptomatic, and~2 means symptomatic.
 From now on, we call each element~$(x, i) \in \Delta$ a site, or a pair cluster-level.
 Using an idea of Harris~\cite{harris_1972} adapted to our model in~\cite[Section~3]{belatreche_belhadji_lanchier_2026}, the process starting from any initial configuration can be constructed from the following collections of independent Poisson processes:
\begin{itemize}
\item {\bf External infections}.
 For~$j = 1, 2$, and all~$(x, i) \in \Delta$ and~$y \sim x$, place an exponential clock with rate~$\lambda_j$ along the edge~$(x, i) \to (y, 1)$, and draw
 $$ \hbox{a type~$j$ arrow} \ (x, i, t) \arrow{$j$} (y, 1, t) \ \hbox{at the times~$t$ the clock rings} $$
 to indicate that if the tail of the arrow is occupied by a~$j$ and each of the levels of~$y$ are occupied by a~0, then the head of the arrow becomes occupied by a~1. \vspace*{4pt}
\item {\bf Internal infections}.
 Similarly, for~$j = 1, 2$, and all~$(x, i) \in \Delta$, place an exponential clock with rate~$\phi_j$ at site~$(x, i)$, and put
 $$ \hbox{a type~$j$ bullet \bul{$j$} at} \ (x, i, t) \ \hbox{at the times~$t$ the clock rings} $$
 to indicate that if the bullet is occupied by a~$j$ and there is at least one~0 in cluster~$x$, then one of the~0s, say the lowest~0 to fix the ideas, turns into a~1. \vspace*{4pt}
\item {\bf Mutations}.
 Place an exponential clock with rate~$\gamma$ at each~$(x, i) \in \Delta$, and put
 $$ \hbox{a circle $\circ$ at} \ (x, i, t) \ \hbox{at the times~$t$ the clock rings} $$
 to indicate that if the circle is occupied by a~1, then this~1 turns into a~2. \vspace*{4pt}
\item {\bf Recoveries}.
 Place an exponential clock with rate one at each~$(x, i) \in \Delta$, and put
 $$ \hbox{a cross $\times$ at} \ (x, i, t) \ \hbox{at the times~$t$ the clock rings} $$
 to indicate that if the cross is occupied by a~1 or a~2, then it turns into a~0.
\end{itemize}


\section{Proof of Theorem~\ref{th:gamma1}}
\begin{figure}[t!]
\centering
\scalebox{0.56}{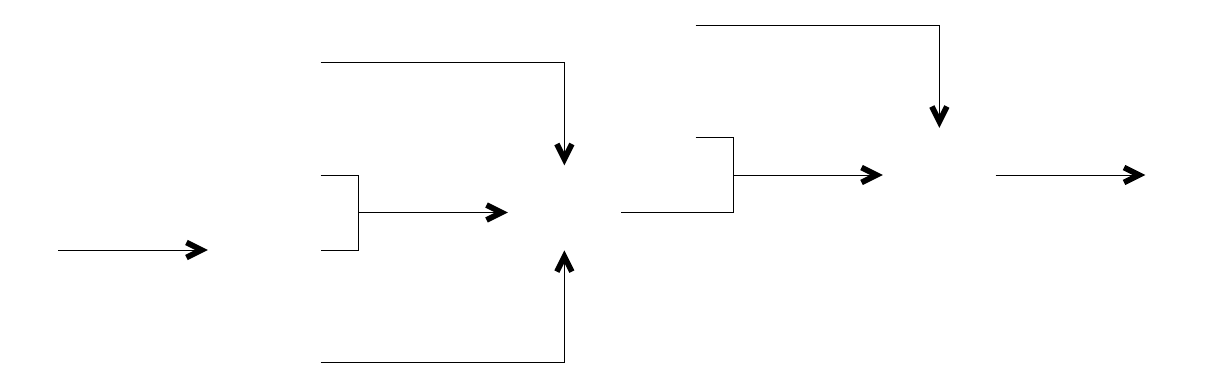}
\caption{\upshape{
 Structure of the proof.
 In the picture, $A \longrightarrow B$ means that~$A$ implies~$B$, while~$A \overset{n}{\longrightarrow} B$ means that there is a good event~$G$ that only depends on the graphical representation in~$\B$, with probability close to one when~$N$ is large, and such that~$A \cap G$ implies~$B$, where~$n$ refers to the lemma where this property is proved.}}
\label{fig:proof1}
\end{figure}

 Theorem~\ref{th:gamma1}.a, which states extinction when the mutation rate is either very small or very large, follows directly from results in~\cite{belatreche_belhadji_lanchier_2026}.
 Indeed,
\begin{itemize}
\item
 Because~$\lambda_1 < \lambda_c$ and~$\phi_1 = 0$, for all values of~$\lambda_2, \phi_2$ and~$N$,~\cite[Theorem~5.b]{belatreche_belhadji_lanchier_2026} implies that there exists~$\gamma_- > 0$ such that the process dies out for all~$\gamma < \gamma_-$. \vspace*{4pt}
\item
 Because~$\lambda_2 = 0$, for all values of~$\lambda_1, \phi_1, \phi_2$ and~$N$,~\cite[Theorem~3.c]{belatreche_belhadji_lanchier_2026} implies that there exists~$\gamma_+ < \infty$ such that the process dies out for all~$\gamma > \gamma_+$.
\end{itemize}
 In fact, these two theorems imply more generally extinction of the stochastic process when the infection rates~$\lambda_2$ and~$\phi_1$ are positive but small. \\
\indent
 The proof of Theorem~\ref{th:gamma1}.b, which states survival for intermediate values of~$\gamma$, and so confirms the presence of at least two phase transitions and the lack of monotonicity with respect to the mutation rate, is more complicated.
 We first define the space-time block
 $$ \B = \{-e_1, 0, e_1 \} \times \{1, 2, \ldots, N \} \times [0, 2T] \quad \hbox{where} \quad T = \ln (\ln (N)), $$
 and where~$e_1$ refers to the first unit vector in~$\Z^d$.
 The goal is to prove that, if there are initially many asymptomatic individuals in cluster~0, then there will be many asymptomatic individuals in cluster~$e_1$ at time~$2T$ with probability close to one and regardless of the graphical representation and configuration of the process outside~$\B$.
 Due to the lack of attractiveness of the process, there is no ``worse case scenario'' for the initial configuration in cluster~$e_1$, so we need to deal with all the possible cases.
 More precisely, we let
 $$ \begin{array}{rclcl}
      H (x, t) & \n = \n & \sum_i \ind \{\xi_t (x, i) = 0 \} & \n = \n & \# \,\hbox{healthy individuals in cluster~$x$}, \vspace*{4pt} \\
      A (x, t) & \n = \n & \sum_i \ind \{\xi_t (x, i) = 1 \} & \n = \n & \# \,\hbox{asymptomatic individuals in cluster~$x$}, \vspace*{4pt} \\
      S (x, t) & \n = \n & \sum_i \ind \{\xi_t (x, i) = 2 \} & \n = \n & \# \,\hbox{symptomatic individuals in cluster~$x$}, \end{array} $$
 and introduce the following ten events:
 $$ \begin{array}{lclrcl}
      A_- & \n = \n & \{A (e_1, 0) < N^{0.6} \},           & \qquad E_0 & \n = \n & \{A (0, 0) \geq N^{0.4} \}, \vspace*{4pt} \\
      A_+ & \n = \n & \{A (e_1, 0) \geq N^{0.6} \},        & \qquad E_1 & \n = \n & \{\forall \,t \leq T, A (0, t) \geq N^{0.2} \}, \vspace*{4pt} \\
      S_0 & \n = \n & \{S (e_1, 0) = 0 \},                 & \qquad E_2 & \n = \n & \{\exists \,t \leq T : 0 < S (e_1, t) < N - N^{0.8} \}, \vspace*{4pt} \\
      S_- & \n = \n & \{0 < S (e_1, 0) < N - N^{0.8} \},   & \qquad E_3 & \n = \n & \{\exists \,t \leq 2T : A (e_1, t) \geq N^{0.6} \}, \vspace*{4pt} \\
      S_+ & \n = \n & \{S (e_1, 0) \geq N - N^{0.8} \}.    & \qquad E_4 & \n = \n & \{A (e_1, 2T) \geq N^{0.4} \}. \end{array} $$
 The events on the left are the possible initial configurations in cluster~$e_1$.
 To prove the theorem, we will show that the following implications occur with probability close to one.
 Starting with many~1s in cluster~0~(event~$E_0$), the number of~1s will remain large until time~$T$~(event~$E_1$), so if cluster~$e_1$ is initially healthy, it will get infected during this time interval and contain at least one~2 by time~$T$~(event~$E_2$).
 If on the contrary cluster~$e_1$ starts with many~2s, some of these~2s will die leaving some space in the cluster for the~1s to spread~(event~$E_2$).
 In any case, there will be again many~1s in cluster~$e_1$ by time~$2T$~(event~$E_3$), and the number of~1s will remain large until that time~(event~$E_4$).
 The structure of the proof is summarized in Figure~\ref{fig:proof1}.
 Before proving the properties listed in the figure, we start with a preliminary result used repeatedly later.
\begin{lemma}
\label{lem:lde}
 Let~$X = \binomial (K, p)$. Then,
 $$ P (X \leq M) \leq \exp (M - Kp (1 - e^{-1})) \quad \hbox{for all} \quad M > 0. $$
\end{lemma}
\begin{proof}
 Writing~$X = X_1 + \cdots + X_K$ as the sum of~$K$ independent Bernoulli random variables with parameter~$p$, the moment generating function is given by
 $$ \begin{array}{rcl}
      M_X (s) & \n = \n & E (e^{sX}) = E (e^{s (X_1 + X_2 + \cdots + X_K)}) = E (e^{sX_1}) E (e^{sX_2}) \cdots E (e^{sX_K}) \vspace*{4pt} \\
              & \n = \n & (P (X_i = 0) + e^s P (X_i = 1))^K = (1 - p (1 - e^s))^K \vspace*{4pt} \\
              & \n \leq \n & \exp (- Kp (1 - e^s)). \end{array} $$
 Then, taking~$s = -1$ and using Markov's inequality, we get
 $$ \begin{array}{rcl}
      P (X \leq M) & \n = \n & P (e^{-X} \geq e^{-M}) \leq e^M E (e^{-X}) \vspace*{4pt} \\
                   & \n = \n & e^M M_X (-1) \leq \exp (M - Kp (1 - e^{-1})) \end{array} $$
 for all~$M > 0$.
 This completes the proof.
\end{proof} \vspace*{8pt} \\
 We can now prove the five lemmas listed in Figure~\ref{fig:proof1}.
\begin{lemma}[Keep some~1s alive at~0]
\label{lem1:keep-0}
 For all~$\ep > 0$, there exists an event~$G_1$ measurable with respect to the graphical representation in~$\B$ such that
 $$ E_0 \cap G_1 \subset E_1 \quad \hbox{and} \quad P (G_1) \geq 1 - \ep/10 \ \hbox{for all~$N$ large}. $$
\end{lemma}
\begin{proof}
 Note that a~1 at a given site can only be removed due to a mutation mark~$\circ$ or a death mark~$\times$ occurring together at rate~$\gamma + 1$.
 In particular, on the event~$E_0$, the number of~1s initially present at cluster~0 that remain~1 until at least time~$T$ dominates
 $$ X = \binomial (N^{0.4}, p) \quad \hbox{where} \quad p = e^{- (\gamma + 1) T} = e^{- (\gamma + 1) \ln (\ln (N))} = (\ln (N))^{- (\gamma + 1)}. $$
 In particular, applying Lemma~\ref{lem:lde} with~$M = N^{0.2}$, for all~$N$ large, we get
 $$ \begin{array}{rcl}
      P (E_1^c \,| \,E_0) & \n = \n &
      P (\exists \,t \leq T : A (0, t) < N^{0.2} \,| \,E_0) \vspace*{4pt} \\ & \n \leq \n &
      P (X \leq N^{0.2}) \leq
     \exp (N^{0.2} - N^{0.4} (1 - e^{-1})(\ln (N))^{- (\gamma + 1)}) \leq \ep/10. \end{array} $$
 Because these estimates only depend on the graphical representation in the block~$\B$, this implies the existence of a good event~$G_1$ that satisfies the lemma.
\end{proof}
\begin{lemma}[Send a~2 onto~$e_1$]
\label{lem1:send-e}
 For all~$\ep > 0$, there exists an event~$G_2$ measurable with respect to the graphical representation in~$\B$ such that
 $$ E_1 \cap S_0 \cap G_2 \subset E_2 \quad \hbox{and} \quad P (G_2) \geq 1 - \ep/10 \ \hbox{for all~$N$ large}. $$
\end{lemma}
\begin{proof}
 To begin with, we denote respectively by
 $$ H = \int_0^T \ind \{H (e_1, s) = N \} \,ds \quad \hbox{and} \quad Y = \poisson (\gamma T/2) $$
 the amount of time cluster~$e_1$ is fully healthy by time~$T$ and the Poisson random variable with parameter~$\gamma T/2$.
 On the event~$E_1$, cluster~$e_1$ gets infected at rate at least~$\lambda_1 N^{0.2}$, while infected individuals recover at rate one, therefore~$E (H)$ is at most~$T / (1 + \lambda_1 N^{0.2})$.
 In addition, each~1 turns into a~2 at rate~$\gamma$.
 In particular, conditioning on whether~$H$ is larger or smaller than~$T/2$, and using~Markov's inequality, we deduce that, for all~$N$ large,
 $$ \begin{array}{rcl}
      P (E_2^c \,| \,E_1 \cap S_0) & \n = \n &
      P (\forall \,t \leq T, S (e_1, t) = 0 \,| \,E_1 \cap S_0) \vspace*{4pt} \\ & \n \leq \n &
      P (H \geq T/2) + P (Y = 0) \leq
      2 E (H) / T + e^{- \gamma T/2} \vspace*{4pt} \\ & \n \leq \n &
      2 / (1 + \lambda_1 N^{0.2}) + (\ln (N))^{- \gamma / 2} \leq \ep/10. \end{array} $$
 As previously, because these estimates only depend on the graphical representation of the process in the space-time block~$\B$, the lemma follows.
\end{proof}
\begin{lemma}[Kill some~2s at~$e_1$]
\label{lem1:kill-e}
 For all~$\ep > 0$, there exists an event~$G_3$ measurable with respect to the graphical representation in~$\B$ such that
 $$ S_+ \cap G_3 \subset E_2 \quad \hbox{and} \quad P (G_3) \geq 1 - \ep/10 \ \hbox{for all~$N$ large}. $$
\end{lemma}
\begin{proof}
 Note that a site initially occupied by a~2 at which there is at least one death mark~$\times$ but no mutation mark~$\circ$ between time~0 and time~1 is not occupied by a~2 at time~1.
 In particular, on the event~$S_+$, the number of~2s at cluster~$e_1$ replaced by a~0 or a~1 at time~1 dominates
 $$ X = \binomial (N - N^{0.8}, p) \quad \hbox{where} \quad p = (1 - e^{-1}) \,e^{- \gamma}. $$
 Applying again Lemma~\ref{lem:lde} but with~$M = N^{0.8}$, we get
 $$ \begin{array}{rcl}
      P (E_2^c \,| \,S_+) & \n \leq \n &
      P (S (e_1, 1) \geq N - N^{0.8} \,| \,S_+) \vspace*{4pt} \\ & \n \leq \n &
      P (X \leq N^{0.8}) \leq
     \exp (N^{0.8} - (N - N^{0.8}) (1 - e^{-1})^2 e^{- \gamma}) \leq \ep/10 \end{array} $$
 for all~$N$ large.
 This completes the proof.
\end{proof}
\begin{lemma}[Produce many~1s at~$e_1$]
\label{lem1:produce-e}
 For all~$\ep > 0$, there exists an event~$G_4$ measurable with respect to the graphical representation in~$\B$ such that
 $$ E_2 \cap A_- \cap G_4 \subset E_3 \quad \hbox{and} \quad P (G_4) \geq 1 - \ep/10 \ \hbox{for all~$N$ large}. $$
\end{lemma}
\begin{proof}
 On the event~$E_2$, with probability one,
 $$ \tau = \inf \{t : 0 < S (e_1, t) < N - N^{0.8} \} < T \quad \hbox{and} \quad \sigma = \tau + N^{- 0.2} < 2T. $$
 If there are more than~$N^{0.6}$ asymptomatic individuals at~$e_1$ at time~$\tau$, we are done.
 Otherwise, we will show that, at time~$\sigma$, the number of~1s at~$e_1$ exceeds~$N^{0.6}$, and so the event~$E_3$ occurs, with probability close to one.
 By our assumption and the definition of time~$\tau$,
\begin{itemize}
\item $S (e_1, \tau) > 0 \ \Longrightarrow \ \xi_{\tau} (e_1, i_0) = 2 \ \hbox{for some level} \ i_0 \in \{1, 2, \ldots, N \}$, \vspace*{4pt}
\item $S (e_1, \tau) < N - N^{0.8}$ and~$A(e_1, \tau) < N^{0.6} \ \Longrightarrow \ H (e_1, \tau) \geq N^{0.8} - N^{0.6} \geq N^{0.7}$.
\end{itemize}
 To prove that many infections originate from the~2 at level~$i_0$, we introduce the event
 $$ B_{i_0} = \{\hbox{no~$\times$ and at least~$N^{0.7}$ type~2 $\bullet$ in~$(e_1, i_0) \times [\tau, \sigma]$} \}. $$
 Let~$\rho$ be the time it takes starting from~$\tau$ to have~$N^{0.7}$ type~2~$\bullet$ at site~$(e_1, i_0)$.
 Using that the death marks~$\times$ occur at rate one, and that, because the internal infections~\bul{2} occur at rate~$\phi_2 = N$, time~$\rho$ has expected value~$N^{0.7} / N = N^{- 0.3}$, it follows from Markov's inequality that
\begin{equation}
\label{eq1:produce-e1}
\begin{array}{rcl}
  P (B_{i_0}^c) & \n \leq \n & P (\hbox{at least one~$\times$ in~$(e_1, i_0) \times [\tau, \sigma]$}) + P (\rho \geq \sigma - \tau) \vspace*{4pt} \\
                & \n \leq \n & (1 - e^{- N^{- 0.2}}) + P (\rho \geq N^{- 0.2}) \vspace*{4pt} \\
                & \n \leq \n & N^{- 0.2} + N^{0.2} E (\rho) = N^{- 0.2} + N^{0.2} N^{- 0.3} = N^{- 0.2} + N^{- 0.1}. \end{array}
\end{equation}
 Now, on the event~$B_{i_0}$, there are at least~$N^{0.7}$ potential infections in the time interval~$[\tau, \sigma]$, and because there are at least~$N^{0.7}$ healthy individuals at time~$\tau$, there must be~$N^{0.7}$ successful infections in the time interval indeed turning a~0 into a~1.
 These successful infections occur in a set~$\Lambda$ of at most~$N^{0.7}$ sites.
 Because mutations~$\circ$ and deaths~$\times$ occur together at rate~$\gamma + 1$, the number of mutations and deaths in~$\Lambda \times [\tau, \sigma]$ is dominated by
 $$ Y = \poisson (\mu) \quad \hbox{where} \quad \mu = (\gamma + 1)(\sigma - \tau) N^{0.7} = (\gamma + 1) N^{0.5}. $$
 In particular, using again Markov's inequality shows that the number of mutations and deaths in the space-time region~$\Lambda \times [\tau, \sigma]$ exceeds~$N^{0.6}$ with probability at most
 $$ P (Y \geq N^{0.6}) \leq N^{- 0.6} E (Y) = N^{- 0.6} (\gamma + 1) N^{0.5} = (\gamma + 1) N^{- 0.1}, $$
 while if there are at most~$N^{0.6}$ mutations and deaths, then
\begin{equation}
\label{eq1:produce-e2}
   A (e_1, \sigma) \geq (\# \hbox{1s created} - \# \hbox{1s removed}) \ \hbox{in} \ \Lambda \times [\tau, \sigma] \geq N^{0.7} - N^{0.6} \geq N^{0.6},
\end{equation}
 and so the event~$E_3$ occurs.
 Combining~\eqref{eq1:produce-e1}--\eqref{eq1:produce-e2}, we conclude that
 $$ P (E_3^c \,| \,E_2 \cap A_-) \leq P (B_{i_0}^c) + P (Y \geq N^{0.6}) \leq N^{- 0.2} + (\gamma + 2) N^{- 0.1} \leq \ep/10 $$
 for all~$N$ large.
 This completes the proof.
\end{proof}
\begin{lemma}[Keep some~1s alive at~$e_1$]
\label{lem1:keep-e}
 For all~$\ep > 0$, there exists an event~$G_5$ measurable with respect to the graphical representation in~$\B$ such that
 $$ A_+ \cap G_5 \subset E_3 \cap G_5 \subset E_4 \quad \hbox{and} \quad P (G_5) \geq 1 - \ep/10 \ \hbox{for all~$N$ large}. $$
\end{lemma}
\begin{proof}
 Because~$A_+ \subset E_3$, the first inclusion is obvious.
 To prove the second inclusion, note that, on the event~$E_3$, we can find a set of~$N^{0.6}$ asymptomatic individuals at cluster~$e_1$ at some time~$t \leq 2T$, and if at least~$N^{0.4}$ of these individuals remain asymptomatic for~$2T - t \leq 2T$ units of time, then the event~$E_4$ occurs.
 Now, using the same argument as in the proof of Lemma~\ref{lem1:keep-0}, the number of~1s at time~$t$ that remain~1 until at least time~$2T$ dominates
 $$ X = \binomial (N^{0.6}, p) \quad \hbox{where} \quad p = e^{- 2 (\gamma + 1) T} = e^{- 2 (\gamma + 1) \ln (\ln (N))} = (\ln (N))^{- 2 (\gamma + 1)}. $$
 Then, applying Lemma~\ref{lem:lde} with~$M = N^{0.4}$, we get
 $$ P (E_4^c \,| \,E_3) \leq P (X \leq N^{0.4}) \leq \exp (N^{0.4} - N^{0.6} (1 - e^{-1})(\ln (N))^{- 2 (\gamma + 1)}) \leq \ep/10 $$
 for all~$N$ large.
 This completes the proof.
\end{proof} \vspace*{8pt} \\
 Using the previous lemmas, we can now prove the theorem. \vspace*{8pt} \\
\begin{proofof}{Theorem~\ref{th:gamma1}.b}
 To deduce survival, we now use a block construction and couple the process with a 2-dependent oriented site percolation model.
 Define the lattice
 $$ \lat = \{(m, n) \in \Z \times \Z_+ : m + n \ \hbox{is even} \}, $$
 which we turn into a directed graph~$\vec \lat$ by placing arrows
 $$ \begin{array}{rcl} (m, n) \to (m', n') & \Longleftrightarrow & |m - m'| = 1 \ \ \hbox{and} \ \ n' = n + 1. \end{array} $$
 Oriented site percolation with parameter~$p$ and range~$R$ is obtained by assuming that sites are open with probability~$p$ and closed with probability~$1 - p$, and that sites that are more than graph distance~$R$ apart are independent.
 Using a contour argument~\cite[Section~10]{durrett_1984} shows that
 $$ p > 1 - 3^{-4(2R + 1)^2} \quad \Longrightarrow \quad P (|C_0| = \infty) > 0, $$
 where~$C_0$ refers to the cluster of sites that can be reached from the origin by a directed path of open sites.
 Returning to our process, we let
 $$ \B_{m, n} = (me_1, 2nT) + \B \qquad \forall \,(m, n) \in \lat, $$
 take $\ep < 3^{- 4 (2R + 1)^2}$ with $R = 2$ in order to have
 $$ \B_{m, n} \cap \B_{m', n'} = \varnothing \quad \hbox{whenever} \quad |m - m'| \vee |n - n'| > R, $$
 and say that~$(m, n) \in \lat$ is a good site if the event
 $$ \begin{array}{rcl}
      E_{m, n} & \n = \n & \{\hbox{cluster~$me_1$ has at least~$N^{0.4}$ asymptomatic individuals at time~$2nT$} \} \vspace*{4pt} \\
               & \n = \n & \{A (me_1, 2nT) \geq N^{0.4} \} \end{array} $$
 occurs.
 Note that~$E_{m, n}$ is the space-time translation of~$E_0$ by vector~$(me_1, 2nT)$.
 Now, define the good event~$G = G_1 \cap \cdots \cap G_5$.
 Then,~$G$ is measurable with respect to the graphical representation in the block space-time~$\B$.
 In addition, it follows from Lemmas~\ref{lem1:keep-0}--\ref{lem1:keep-e} that
\begin{equation}
\label{eq:th11}
\begin{array}{rcll}
  E_0 \cap G & \n \subset \n &  E_1 \cap G                                               & \quad \hbox{by Lemma~\ref{lem1:keep-0}} \vspace*{4pt} \\
             & \n \subset \n & (E_1 \cap S_0 \cap G) \cup (S_- \cap G) \cup (S_+ \cap G) & \quad \hbox{by $\{S_0, S_-, S_+ \}$ partition} \qquad \vspace*{4pt} \\
             & \n \subset \n & (E_2 \cap G) \cup (S_- \cap G) \cup (S_+ \cap G)          & \quad \hbox{by Lemma~\ref{lem1:send-e}} \vspace*{4pt} \\
             & \n \subset \n &  E_2 \cap G                                               & \quad \hbox{by Lemma~\ref{lem1:kill-e}} \vspace*{4pt} \\
             & \n \subset \n & (E_2 \cap A_- \cap G) \cup (A_+ \cap G)                   & \quad \hbox{by $\{A_-, A_+ \}$ partition} \vspace*{4pt} \\
             & \n \subset \n & (E_3 \cap G) \cup (A_+ \cap G)                            & \quad \hbox{by Lemma~\ref{lem1:produce-e}} \vspace*{4pt} \\
             & \n \subset \n &  E_4 = E_{1, 1}                                           & \quad \hbox{by Lemma~\ref{lem1:keep-e}}. \end{array}
\end{equation}
 Taking the complement, using the union bound, and applying the second part of Lemmas~\ref{lem1:keep-0}--\ref{lem1:keep-e} also imply that, for all~$N$ sufficiently large,
\begin{equation}
\label{eq:th12}
  P (G) = 1 - P \bigg(\bigcup_{i = 1}^5 G_i^c \bigg) \geq 1 - \sum_{i = 1}^5 P (G_i^c) \geq 1 - 5 (\ep/10) = 1 - \ep/2.
\end{equation}
 Finally, using symmetry and translation invariance in space and time of the graphical representation, and combining~\eqref{eq:th11}--\eqref{eq:th12}, we deduce the existence of a collection of good events~$G_{m, n}$, measurable with respect to the graphical representation in the block~$\B_{m, n}$, such that
 $$ E_{m, n} \cap G_{m, n} \subset E_{m - 1, n + 1} \cap E_{m + 1, n + 1} \quad \hbox{and} \quad P (G_{m, n}) \geq 1 - \ep \ \hbox{for all~$N$ large}. $$
 This and~\cite[Theorem~4.3]{durrett_1995} imply that the good sites dominate the open sites in a~2-dependent oriented site percolation process with parameter $1 - \epsilon$.
 In particular, by our choice of~$\ep$, the open sites percolate, so the good sites percolate as well, which proves survival.
\end{proofof}


\section{Proof of Theorem~\ref{th:gamma2}}
\begin{figure}[t!]
\centering
\scalebox{0.56}{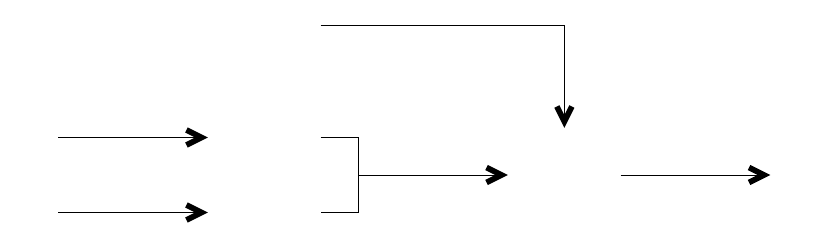}
\caption{\upshape{
 Structure of the proof.
 In the picture, $A \longrightarrow B$ means that~$A$ implies~$B$, while~$A \overset{n}{\longrightarrow} B$ means that there is a good event~$G$ that only depends on the graphical representation in~$\B$, with probability close to one when~$N$ is large, and such that~$A \cap G$ implies~$B$, where~$n$ refers to the lemma where this property is proved.}}
\label{fig:proof2}
\end{figure}

 As previously, the first part of the theorem follows directly from~\cite{belatreche_belhadji_lanchier_2026}.
\begin{itemize}
\item
 Because~$\lambda_1 = 0$, for all values of~$\lambda_2, \phi_1, \phi_2$ and~$N$,~\cite[Theorem~3.b]{belatreche_belhadji_lanchier_2026} implies that there exists~$\gamma_- > 0$ such that the process dies out for all~$\gamma < \gamma_-$. \vspace*{4pt}
\item
 Because~$\lambda_2 < \lambda_c$ and~$\phi_2 = 0$, for all values of~$\lambda_1, \phi_1$ and~$N$,~\cite[Theorem~5.c]{belatreche_belhadji_lanchier_2026} implies that there exists~$\gamma_+ < \infty$ such that the process dies out for all~$\gamma > \gamma_+$.
\end{itemize}
 These two theorems again imply more generally extinction of the stochastic process when the infection rates~$\lambda_1$ and~$\phi_2$ are positive but small. \\
\indent
 Theorem~\ref{th:gamma2}.b, which states survival for intermediate values of the mutation rate and the presence of at least two phase transitions, but for other values of the infection rates, again relies on a block construction, though the general strategy strongly differs from Theorem~\ref{th:gamma1}.b.
 More precisely, the space-time block~$\B$ is now defined using a much larger time scale:
 $$ \B = \{-e_1, 0, e_1 \} \times \{1, 2, \ldots, N \} \times [0, T] \quad \hbox{where} \quad T = N^{0.4}, $$
 and the goal is to prove that, if we start with at least one asymptomatic individual in cluster~0, then there will be at least one asymptomatic individual in cluster~$e_1$ at time~$T$ with probability close to one when the cluster size~$N$ is large, regardless of the graphical representation and configuration of the process outside~$\B$.
 As previously, there is no ``worst-case scenario'' for the initial configuration in cluster~$e_1$ so we must deal with all the possible cases.
 Define the seven events
 $$ \begin{array}{lclrcl}
      A_0 & \n = \n & \{A (e_1, 0) = 0 \},           & \qquad E_0 & \n = \n & \{A (0, 0) \neq 0 \}, \vspace*{4pt} \\
      A_+ & \n = \n & \{A (e_1, 0) \neq 0 \},        & \qquad E_1 & \n = \n & \{\forall \,t \leq T, A (0, t) \neq 0 \}, \vspace*{4pt} \\
                                                   &&& \qquad E_2 & \n = \n & \{\exists \, t \leq \sqrt{T} : H (e_1, t) = N \}, \vspace*{4pt} \\
                                                   &&& \qquad E_3 & \n = \n & \{\exists \,t \leq T : A (e_1, t) \neq 0 \}, \vspace*{4pt} \\
                                                   &&& \qquad E_4 & \n = \n & \{A (e_1, T) \neq 0 \}. \end{array} $$
 To prove the theorem, we will first show that, with probability close to one when the cluster sie is large, starting with at least one~1 in cluster~0~(event~$E_0$), there will be at least one~1 in the cluster at all times until time~$T$~(event~$E_1$).
 In addition, if cluster~$e_1$ starts with no~1s, it will become healthy by time~$\sqrt{T}$~(event~$E_2$).
 Then, either cluster~$e_1$ starts with at least one~1~(event~$E_3$), or it starts with no~1s, in which case many~1s in cluster~0 turn into a~2, and one of these~2s will infect cluster~$e_1$ between time~$\sqrt{T}$ and time~$T$~(event~$E_3$).
 Once we have a~1 in cluster~$e_1$, there will be at least one~1 in the cluster at all times until time~$T$~(event~$E_4$).
 The structure of the proof is summarized in Figure~\ref{fig:proof2}, and consists of the following four lemmas.
\begin{lemma}[Keep a~1 alive at~0]
\label{lem2:keep-0}
 For all~$\ep > 0$, there exists an event~$G_1$ measurable with respect to the graphical representation in~$\B$ such that
 $$ E_0 \cap G_1 \subset E_1 \quad \hbox{and} \quad P (G_1) \geq 1 - \ep/8 \ \hbox{for all~$N$ large}. $$
\end{lemma}
\begin{proof}
 The idea is to define a set-valued process~$\Xi_t$ that keeps track of a subset of asymptomatic individuals as follows.
 The process starts at
 $$ \Xi_0 = \{i_1 \} \quad \hbox{where} \quad i_1 = \min \{i : \xi_0 (0, i) = 1 \}. $$
 Note that such a level~$i_1$ exists because~$A (0, 0) \neq 0$.
 Then, define
 $$ \begin{array}{rcl}
    \sigma_1 & \n = \n & \hbox{time of the first death~$\times$ in the cluster at a level~$\neq i_0$}, \vspace*{4pt} \\
    \sigma_2 & \n = \n & \hbox{time of the first internal infection~\bul{1} after time~$\sigma_1$ at site~$(0, i_1)$}, \vspace*{4pt} \\
    \sigma_3 & \n = \n & \hbox{time of the first mutation~$\circ$ or death mark~$\times$ at site~$(0, i_1)$}. \end{array} $$
 On the event that~$\sigma_3 < \sigma_2$, we define
 $$ \Xi_t = \{i_1 \} \ \ \forall \,t \in [0, \sigma_3), \qquad \Xi_t = \varnothing \ \ \forall \,t \in [\sigma_3, \infty), $$
 to indicate that there might be no~1s left in the cluster starting from time~$\sigma_3$.
 In contrast, on the event that~$\sigma_3 > \sigma_2$, there is at least one healthy level at time~$\sigma_1$ and the~1 at level~$i_1$ attempts to infect a healthy level at time~$\sigma_2 > \sigma_1$.
 This forces the existence of a healthy level, say~$j_1$, becoming infected at some time~$\sigma_4$ in the time interval~$(\sigma_1, \sigma_2]$, and we put
 $$ \sigma_5 = \hbox{time of the first~$\circ$ or~$\times$ at level~$i_1$ or~$j_1$ after time~$\sigma_4$}. $$
 Letting~$i_2 \in \{i_1, j_1 \}$ be the level where the~$\circ$ or~$\times$ occurs, we define
 $$ \Xi_t = \{i_1 \} \ \ \forall \,t \in [0, \sigma_4), \qquad \Xi_t = \{i_1, j_1 \} \ \ \forall \,t \in [\sigma_4, \sigma_5), \qquad \Xi_{\sigma_5} = \{i_2 \}. $$
 Then, the process is defined by reiterating the same procedure but starting from~$\{i_2 \}$ at time~$\sigma_5$ rather than from~$\{i_1 \}$ at time~0, and so on.
 This continues until the asymptomatic individual at level~$i_K$ mutates or dies before trying to infect another site, at which time the process gets trapped in state~$\varnothing$.
 By construction, each level in~$\Xi_t$ is in state~1 at time~$t$, therefore
\begin{equation}
\label{eq2:keep-01}
  A (0, t) \geq \card (\Xi_t) \quad \hbox{and} \quad \tau = \inf \{t : A (0, t) = 0 \} \geq \inf \{t : \Xi_t = \varnothing \} = \sigma.
\end{equation}
 Now, by the memoryless and superposition properties,
 $$ P (\sigma_3 > \sigma_2) = P (\sigma_3 > \sigma_2 \,| \,\sigma_3 > \sigma_1) P (\sigma_3 > \sigma_1) = \bigg(\frac{N}{N + \gamma + 1} \bigg) \bigg(\frac{N - 1}{N + \gamma} \bigg) \geq 1 - N^{- 0.8} $$
 for all~$N$ large, so the number of steps~$K$ until the process hits~$\varnothing$ dominates stochastically the geometric random variable~$Y = \geometric (N^{- 0.8})$, and we have
\begin{equation}
\label{eq2:keep-02}
\begin{array}{rcl}
  P (K \leq N^{0.6}) & \n \leq \n &
  P (Y \leq N^{0.6}) = 1 - P (Y > N^{0.6}) \vspace*{4pt} \\ & \n = \n &
  1 - (1 - N^{- 0.8})^{N^{0.6}} \leq 1 - \exp (- 2N^{- 0.2}) \leq 2N^{- 0.2}. \end{array}
\end{equation}
 Also, because~$\sigma_5 \geq \sigma_5 - \sigma_4 = \exponential (2 (\gamma + 1))$, the time to extinction~$\sigma$ dominates
 $$ Z_1 + Z_2 + \cdots + Z_K \quad \hbox{with} \quad Z_i = \exponential (2 (\gamma + 1)) \ \hbox{independent}, $$
 from which it follows that,
 for all~$s > 0$,
 $$ E (e^{-s \sigma} \,| \,K > n) \leq E (e^{- s (Z_1 + Z_2 + \cdots + Z_n)}) = \prod_{i = 1}^n \,E (e^{-s Z_i}) = \bigg(1 + \frac{s}{2 (\gamma + 1)} \bigg)^{-n}. $$
 In particular, taking~$s = 2 (\gamma + 1)(e - 1) > 0$, we get
\begin{equation}
\label{eq2:keep-03}
  E (e^{- 2 (\gamma + 1)(e - 1) \sigma} \,| \,K > n) \leq \bigg(1 + \frac{2 (\gamma + 1)(e - 1)}{2 (\gamma + 1)} \bigg)^{-n} = e^{-n}.
\end{equation}
 Using~\eqref{eq2:keep-01}--\eqref{eq2:keep-03} and applying Markov's inequality, we conclude that
 $$ \begin{array}{rcl}
      P (E_2^c \,| \,E_1) & \n = \n &
      P (\tau \leq N^{0.4}) \leq
      P (\sigma \leq N^{0.4}) \vspace*{4pt} \\ & \n \leq \n &
      P (K \leq N^{0.6}) + P (\sigma \leq N^{0.4} \,| \,K > N^{0.6}) \vspace*{4pt} \\ & \n \leq \n &
      2 N^{- 0.2} + P (e^{- 2 (\gamma + 1)(e - 1) \sigma} \geq e^{- 2 (\gamma + 1)(e - 1) N^{0.4}} \,| \,K > N^{0.6}) \vspace*{4pt} \\ & \n \leq \n &
      2 N^{- 0.2} + e^{N^{0.4}} E (e^{- 2 (\gamma + 1)(e - 1) \sigma} \,| \,K > N^{0.6}) \vspace*{4pt} \\ & \n \leq \n &
      2 N^{- 0.2} + \exp (N^{0.4} - N^{0.6}) \leq \ep/8 \end{array} $$
 for all~$N$ large.
 This completes the proof.
\end{proof}
\begin{lemma}[Kill all the~2s at~$e_1$]
\label{lem2:kill-e}
 For all~$\ep > 0$, there exists an event~$G_2$ measurable with respect to the graphical representation in~$\B$ such that
 $$ A_0 \cap G_2 \subset E_2 \quad \hbox{and} \quad P (G_2) \geq 1 - \ep/8 \ \hbox{for all~$N$ large}. $$
\end{lemma}
\begin{proof}
 The idea is to study the time it takes to have at least one death mark~$\times$ at each level of the cluster.
 To this end, we let~$\tau_0 = 0$ and
 $$ \tau_1 < \tau_2 < \cdots < \tau_N \ \ \hbox{be the times of the first~$\times$ at the~$N$ levels of cluster~$e_1$}, $$
 ordered from the smallest to the largest.
 Notice that the cluster becomes healthy by time~$\tau_N$ in general, and exactly at time~$\tau_N$ when starting from~$N$ symptomatic individuals.
 Now, by the memoryless and superposition properties for the exponential distribution,
 $$ E (\tau_N) = \sum_{i = 0}^{N - 1} E (\tau_{i + 1} - \tau_i) = \frac{1}{1} + \frac{1}{2} + \frac{1}{3} + \cdots + \frac{1}{N} \leq 1 + \ln (N). $$
 In particular, using Markov's inequality, we conclude that
 $$ P (E_2^c \,| \,A_0) \leq P (\tau_N \geq N^{0.2}) \leq N^{- 0.2} E (\tau_N) \leq N^{- 0.2} (1 + \ln (N)) \leq \ep/8 $$
 for all~$N$ large.
 This completes the proof.
\end{proof}
\begin{lemma}[Send a~1 onto~$e_1$]
\label{lem2:send-e}
 For all~$\ep > 0$, there exists an event~$G_3$ measurable with respect to the graphical representation in~$\B$ such that
 $$ E_1 \cap E_2 \cap G_3 \subset E_3 \quad \hbox{and} \quad P (G_3) \geq 1 - \ep/8 \ \hbox{for all~$N$ large}. $$
\end{lemma}
\begin{proof}
 Because, on the event~$E_1$, cluster~0 is occupied by at least one~1 until time~$T$, the number of times a~1 turns into a~2 between time~$\sqrt{T}$ and time~$T$ dominates
 $$ Y_1 = \poisson (\gamma (T - \sqrt{T})). $$
 Moreover, each such~2 attempts to infect cluster~$e_1$ at rate~$\lambda_2$ and recovers at rate one, so the probability of at least one attempt is~$\lambda_2 / (\lambda_2 + 1)$.
 In particular, by the thinning property, the number of attempts to infect cluster~$e_1$ during the time window~$[\sqrt{T}, T]$ dominates
 $$ Y_2 = \poisson \bigg(\frac{\lambda_2}{\lambda_2 + 1} \ \gamma (T - \sqrt{T}) \bigg) = \poisson \bigg(\frac{\lambda_2}{\lambda_2 + 1} \ \gamma (N^{0.4} - N^{0.2}) \bigg). $$
 Now, on the event~$E_2$ that cluster~$e_1$ is healthy by time~$\sqrt{T} = N^{0.2}$, the cluster remains healthy until time~$T$ only if no~2s attempt to infect the cluster until that time, so
 $$ P (E_3^c \,| \,E_1 \cap E_2) \leq P (Y_2 = 0) = \exp \left (- \frac{\lambda_2}{\lambda_2 + 1} \ \gamma (N^{0.4} - N^{0.2})\right ) \leq \ep/8 $$
 for all~$N$ large.
 This completes the proof.
\end{proof}
\begin{lemma}[Keep a~1 alive at~$e_1$]
\label{lem2:keep-e}
 For all~$\ep > 0$, there exists an event~$G_4$ measurable with respect to the graphical representation in~$\B$ such that
 $$ A_+ \cap G_4 \subset E_3 \cap G_4 \subset E_4 \quad \hbox{and} \quad P (G_4) \geq 1 - \ep/8 \ \hbox{for all~$N$ large}. $$
\end{lemma}
\begin{proof}
 Because~$A_+ \subset E_3$, the first inclusion is obvious, while Lemma~\ref{lem2:keep-0}, together with the translation invariance of the graphical representation, gives the second inclusion.
\end{proof} \vspace*{8pt} \\
 Using the previous lemmas, we can now prove the theorem. \vspace*{8pt} \\
\begin{proofof}{Theorem~\ref{th:gamma2}.b}
 As previously, we couple the process with~2-dependent oriented site percolation on the directed graph~$\vec \lat$.
 We define
 $$ \B_{m, n} = (me_1, nT) + \B \qquad \forall \,(m, n) \in \lat $$
 and take~$\ep < 3^{- 4 (2R + 1)^2}$ with $R = 2$, but now say that~$(m, n) \in \lat$ is a good site if
 $$ \begin{array}{rcl}
      E_{m, n} & \n = \n & \{\hbox{cluster~$me_1$ has at least one asymptomatic individual at time~$nT$} \} \vspace*{4pt} \\
               & \n = \n & \{A (me_1, nT) \neq 0 \} \end{array} $$
 occurs.
 As previously, the event~$G = G_1 \cap \cdots \cap G_4$ is measurable with respect to the graphical representation in the block~$\B$, and it follows from Lemmas~\ref{lem2:keep-0}--\ref{lem2:keep-e} that
\begin{equation}
\label{eq:th21}
\begin{array}{rcll}
  E_0 \cap G & \n \subset \n &  E_1 \cap G                             & \quad \hbox{by Lemma~\ref{lem2:keep-0}} \vspace*{4pt} \\
             & \n \subset \n & (E_0 \cap A_0 \cap G) \cup (A_+ \cap G) & \quad \hbox{by $\{A_0, A_+ \}$ partition} \qquad \vspace*{4pt} \\
             & \n \subset \n & (E_1 \cap E_2 \cap G) \cup (A_+ \cap G) & \quad \hbox{by Lemma~\ref{lem2:kill-e}} \vspace*{4pt} \\
             & \n \subset \n & (E_3 \cap G) \cup (A_+ \cap G)          & \quad \hbox{by Lemma~\ref{lem2:send-e}} \vspace*{4pt} \\
             & \n \subset \n &  E_3 \cap G \subset E_4 = E_{1, 1}      & \quad \hbox{by Lemma \ref{lem2:keep-e}}, \end{array}
\end{equation}
 and that, for all~$N$ sufficiently large,
\begin{equation}
\label{eq:th22}
  P (G) = 1 - P \bigg(\bigcup_{i = 1}^4 G_i^c \bigg) \geq 1 - \sum_{i = 1}^4 P (G_i^c) \geq 1 - 4 (\ep/8) = 1 - \ep/2.
\end{equation}
 Using again symmetry and translation invariance in space and time of the graphical representation, and combining~\eqref{eq:th21}--\eqref{eq:th22}, we deduce the existence of a collection of good events~$G_{m, n}$, measurable with respect to the graphical representation in the block~$\B_{m, n}$, such that
 $$ E_{m, n} \cap G_{m, n} \subset E_{m - 1, n + 1} \cap E_{m + 1, n + 1} \quad \hbox{and} \quad P (G_{m, n}) \geq 1 - \ep \ \hbox{for all~$N$ large}. $$	
 As in the proof of Theorem \ref{th:gamma1}.b, survival follows from~\cite[Theorem~4.3]{durrett_1995}.
\end{proofof}


\end{document}

%% file: gamma.pdf_t
\begin{picture}(0,0)%
\includegraphics{gamma.pdf}%
\end{picture}%
\setlength{\unitlength}{3947sp}%
\begingroup\makeatletter\ifx\SetFigFont\undefined%
\gdef\SetFigFont#1#2#3#4#5{%
  \reset@font\fontsize{#1}{#2pt}%
  \fontfamily{#3}\fontseries{#4}\fontshape{#5}%
  \selectfont}%
\fi\endgroup%
\begin{picture}(11902,9088)(-331,-7934)
\put(1351,989){\makebox(0,0)[b]{\smash{{\SetFigFont{12}{14.4}{\familydefault}{\mddefault}{\updefault}{\color[rgb]{0,0,0}Time 1}%
}}}}
\put(4201,989){\makebox(0,0)[b]{\smash{{\SetFigFont{12}{14.4}{\familydefault}{\mddefault}{\updefault}{\color[rgb]{0,0,0}Time 10}%
}}}}
\put(7051,989){\makebox(0,0)[b]{\smash{{\SetFigFont{12}{14.4}{\familydefault}{\mddefault}{\updefault}{\color[rgb]{0,0,0}Time 100}%
}}}}
\put(9901,989){\makebox(0,0)[b]{\smash{{\SetFigFont{12}{14.4}{\familydefault}{\mddefault}{\updefault}{\color[rgb]{0,0,0}Time 1000}%
}}}}
\put(5551,-7861){\makebox(0,0)[rb]{\smash{{\SetFigFont{12}{14.4}{\familydefault}{\mddefault}{\updefault}{\color[rgb]{0,0,0}Number of infected individuals per cluster:}%
}}}}
\put(-149,-6211){\rotatebox{90.0}{\makebox(0,0)[b]{\smash{{\SetFigFont{12}{14.4}{\familydefault}{\mddefault}{\updefault}{\color[rgb]{0,0,0}$\gamma = 2.00$}%
}}}}}
\put(11401,-6211){\rotatebox{270.0}{\makebox(0,0)[b]{\smash{{\SetFigFont{12}{14.4}{\familydefault}{\mddefault}{\updefault}{\color[rgb]{0,0,0}Theorem \ref{th:gamma1}.b}%
}}}}}
\put(-149,-511){\rotatebox{90.0}{\makebox(0,0)[b]{\smash{{\SetFigFont{12}{14.4}{\familydefault}{\mddefault}{\updefault}{\color[rgb]{0,0,0}$\gamma = 0.80$}%
}}}}}
\put(11401,-511){\rotatebox{270.0}{\makebox(0,0)[b]{\smash{{\SetFigFont{12}{14.4}{\familydefault}{\mddefault}{\updefault}{\color[rgb]{0,0,0}Theorem \ref{th:gamma1}.a}%
}}}}}
\put(-149,-3361){\rotatebox{90.0}{\makebox(0,0)[b]{\smash{{\SetFigFont{12}{14.4}{\familydefault}{\mddefault}{\updefault}{\color[rgb]{0,0,0}$\gamma = 5.00$}%
}}}}}
\put(11401,-3361){\rotatebox{270.0}{\makebox(0,0)[b]{\smash{{\SetFigFont{12}{14.4}{\familydefault}{\mddefault}{\updefault}{\color[rgb]{0,0,0}Theorem \ref{th:gamma1}.a}%
}}}}}
\put(5926,-7861){\makebox(0,0)[lb]{\smash{{\SetFigFont{12}{14.4}{\familydefault}{\mddefault}{\updefault}{\color[rgb]{0,0,0}= 0}%
}}}}
\put(6601,-7861){\makebox(0,0)[lb]{\smash{{\SetFigFont{12}{14.4}{\familydefault}{\mddefault}{\updefault}{\color[rgb]{0,0,0}= 1}%
}}}}
\put(7276,-7861){\makebox(0,0)[lb]{\smash{{\SetFigFont{12}{14.4}{\familydefault}{\mddefault}{\updefault}{\color[rgb]{0,0,0}= 2}%
}}}}
\put(7951,-7861){\makebox(0,0)[lb]{\smash{{\SetFigFont{12}{14.4}{\familydefault}{\mddefault}{\updefault}{\color[rgb]{0,0,0}= 3}%
}}}}
\end{picture}%

%% file: proof1.pdf_t
\begin{picture}(0,0)%
\includegraphics{proof1.pdf}%
\end{picture}%
\setlength{\unitlength}{3947sp}%
\begingroup\makeatletter\ifx\SetFigFont\undefined%
\gdef\SetFigFont#1#2#3#4#5{%
  \reset@font\fontsize{#1}{#2pt}%
  \fontfamily{#3}\fontseries{#4}\fontshape{#5}%
  \selectfont}%
\fi\endgroup%
\begin{picture}(9630,3076)(-2114,-8570)
\put(  1,-8461){\makebox(0,0)[b]{\smash{{\SetFigFont{20}{24.0}{\familydefault}{\mddefault}{\updefault}$S_+$}}}}
\put(1276,-8236){\makebox(0,0)[b]{\smash{{\SetFigFont{14}{16.8}{\familydefault}{\mddefault}{\updefault}\ref{lem1:kill-e}}}}}
\put(  1,-6061){\makebox(0,0)[b]{\smash{{\SetFigFont{20}{24.0}{\familydefault}{\mddefault}{\updefault}$S_-$}}}}
\put(3001,-5761){\makebox(0,0)[b]{\smash{{\SetFigFont{20}{24.0}{\familydefault}{\mddefault}{\updefault}$A_+$}}}}
\put(-2099,-7561){\makebox(0,0)[b]{\smash{{\SetFigFont{20}{24.0}{\familydefault}{\mddefault}{\updefault}$E_0$}}}}
\put(  1,-7561){\makebox(0,0)[b]{\smash{{\SetFigFont{20}{24.0}{\familydefault}{\mddefault}{\updefault}$E_1$}}}}
\put(  1,-6961){\makebox(0,0)[b]{\smash{{\SetFigFont{20}{24.0}{\familydefault}{\mddefault}{\updefault}$S_0$}}}}
\put(2401,-7261){\makebox(0,0)[b]{\smash{{\SetFigFont{20}{24.0}{\familydefault}{\mddefault}{\updefault}$E_2$}}}}
\put(3001,-6661){\makebox(0,0)[b]{\smash{{\SetFigFont{20}{24.0}{\familydefault}{\mddefault}{\updefault}$A_-$}}}}
\put(5401,-6961){\makebox(0,0)[b]{\smash{{\SetFigFont{20}{24.0}{\familydefault}{\mddefault}{\updefault}$E_3$}}}}
\put(7501,-6961){\makebox(0,0)[b]{\smash{{\SetFigFont{20}{24.0}{\familydefault}{\mddefault}{\updefault}$E_4$}}}}
\put(-1124,-7336){\makebox(0,0)[b]{\smash{{\SetFigFont{14}{16.8}{\familydefault}{\mddefault}{\updefault}\ref{lem1:keep-0}}}}}
\put(6376,-6736){\makebox(0,0)[b]{\smash{{\SetFigFont{14}{16.8}{\familydefault}{\mddefault}{\updefault}\ref{lem1:keep-e}}}}}
\put(1276,-7036){\makebox(0,0)[b]{\smash{{\SetFigFont{14}{16.8}{\familydefault}{\mddefault}{\updefault}\ref{lem1:send-e}}}}}
\put(4276,-6736){\makebox(0,0)[b]{\smash{{\SetFigFont{14}{16.8}{\familydefault}{\mddefault}{\updefault}\ref{lem1:produce-e}}}}}
\end{picture}%

%% file: proof2.pdf_t
\begin{picture}(0,0)%
\includegraphics{proof2.pdf}%
\end{picture}%
\setlength{\unitlength}{3947sp}%
\begingroup\makeatletter\ifx\SetFigFont\undefined%
\gdef\SetFigFont#1#2#3#4#5{%
  \reset@font\fontsize{#1}{#2pt}%
  \fontfamily{#3}\fontseries{#4}\fontshape{#5}%
  \selectfont}%
\fi\endgroup%
\begin{picture}(6630,1876)(-2114,-7670)
\put(-2099,-6961){\makebox(0,0)[b]{\smash{{\SetFigFont{20}{24.0}{\familydefault}{\mddefault}{\updefault}$E_0$}}}}
\put(  1,-6961){\makebox(0,0)[b]{\smash{{\SetFigFont{20}{24.0}{\familydefault}{\mddefault}{\updefault}$E_1$}}}}
\put(-2099,-7561){\makebox(0,0)[b]{\smash{{\SetFigFont{20}{24.0}{\familydefault}{\mddefault}{\updefault}$A_0$}}}}
\put(  1,-7561){\makebox(0,0)[b]{\smash{{\SetFigFont{20}{24.0}{\familydefault}{\mddefault}{\updefault}$E_2$}}}}
\put(  1,-6061){\makebox(0,0)[b]{\smash{{\SetFigFont{20}{24.0}{\familydefault}{\mddefault}{\updefault}$A_+$}}}}
\put(-1124,-6736){\makebox(0,0)[b]{\smash{{\SetFigFont{14}{16.8}{\familydefault}{\mddefault}{\updefault}\ref{lem2:keep-0}}}}}
\put(-1124,-7336){\makebox(0,0)[b]{\smash{{\SetFigFont{14}{16.8}{\familydefault}{\mddefault}{\updefault}\ref{lem2:kill-e}}}}}
\put(1276,-7036){\makebox(0,0)[b]{\smash{{\SetFigFont{14}{16.8}{\familydefault}{\mddefault}{\updefault}\ref{lem2:send-e}}}}}
\put(2401,-7261){\makebox(0,0)[b]{\smash{{\SetFigFont{20}{24.0}{\familydefault}{\mddefault}{\updefault}$E_3$}}}}
\put(4501,-7261){\makebox(0,0)[b]{\smash{{\SetFigFont{20}{24.0}{\familydefault}{\mddefault}{\updefault}$E_4$}}}}
\put(3376,-7036){\makebox(0,0)[b]{\smash{{\SetFigFont{14}{16.8}{\familydefault}{\mddefault}{\updefault}\ref{lem2:keep-e}}}}}
\end{picture}%